\documentclass [twoside,reqno,draft,12pt] {amsart}

\usepackage{amsfonts}
\usepackage{amssymb}
\usepackage{a4}

\newtheorem{thm}{Theorem}[section]
\newtheorem{cor}[thm]{Corollary}
\newtheorem{lem}[thm]{Lemma}

\theoremstyle{definition}

\newtheorem{exmp}{Example}[section]

\theoremstyle{remark}
\newtheorem{rem}{Remark}[section]

\numberwithin{equation}{section}

\newtheorem{definition}{Definition}[section]

\newtheorem{remark}[definition]{Remark}

\renewcommand{\Re}{\hbox{Re}\,}
\renewcommand{\Im}{\hbox{Im}\,}

\newcommand{\C}{\mathbb{C}}

\newcommand{\R}{\mathbb{R}}

\newcommand{\pM}{{\partial M}}

\newcommand{\la}{{\lambda}}

\def\hat{\widehat}
\def\tilde{\widetilde}
\def \bfo {\begin {eqnarray*} }
\def \efo {\end {eqnarray*} }
\def \ba {\begin {eqnarray*} }
\def \ea {\end {eqnarray*} }
\def \beq {\begin {eqnarray}}
\def \eeq {\end {eqnarray}}

\def \det {\hbox{det}}
\def\bra{\langle}
\def\cet{\rangle}
\def \e {\varepsilon}
\def \p {\partial}

\def\hat{\widehat}
\def\tilde{\widetilde}
\def \bfo {\begin {eqnarray*} }
\def \efo {\end {eqnarray*} }
\def \ba {\begin {eqnarray*} }
\def \ea {\end {eqnarray*} }
\def \beq {\begin {eqnarray}}
\def \eeq {\end {eqnarray}}

\def \det {\hbox{det}}
\def\bra{\langle}
\def\cet{\rangle}

\def \e {\varepsilon}
\def \p {\partial}

\begin{document}

 \title{Inverse spectral problems on a closed manifold}

%% First author

\author[Krupchyk]{Katsiaryna Krupchyk}

\address
        {K. Krupchyk, Department of Physics and Mathematics\\
        University of Joensuu\\
         P.O. Box 111\\
         FI-80101 Joensuu\\
         Finland}

\email{katya.krupchyk@joensuu.fi}

%% Second author

\author[Kurylev]{Yaroslav Kurylev}

\address
        {Y. Kurylev, Department of Mathematics\\
University College London\\
 Gower Street\\
London\\
WC1E 5BT\\
 UK}

\email{Y.Kurylev@ucl.ac.uk}

%% Third author

\author[Lassas]{Matti Lassas}

\address
        {M. Lassas, Institute of Mathematics
\\
P.O.Box 1100\\
 02015 Helsinki
 University of Technology\\
         Finland}

\email{Matti.Lassas@hut.fi}

%\date{Draft version 31.8.2007}

\maketitle

 \begin{abstract}
 In this paper we consider two inverse   problems on a closed connected
Riemannian manifold $(M,g)$. The first one is a direct analog
of the Gel'fand inverse boundary spectral problem. To formulate it,
 assume that
$M$ is divided by a hypersurface $\Sigma$ into two components
and we know the eigenvalues $\la_j$ of the Laplace operator on $(M,g)$
and also the Cauchy data, on $\Sigma$, of the corresponding eigenfunctions
$\phi_j$, i.e. $\phi_j|_{\Sigma},\,\p_\nu \phi_j|_{\Sigma}$, where $\nu$
is the normal to $\Sigma$. We prove that these data determine $(M,g)$
uniquely, i.e. up to an isometry. In the second problem we are given
much less data, namely, $\la_j$ and  $\phi_j|_{\Sigma}$ only.
However, if $\Sigma$ consists of at least two components,
$\Sigma_1, \Sigma_2$, we are still
able to determine $(M,g)$ assuming some conditions on $M$ and $\Sigma$.
These conditions are formulated in terms of the spectra
of the manifolds with boundary obtained by cutting $M$ along $\Sigma_i,
\, i=1,2,$ and are of a generic nature. We consider also some other
inverse problems on $M$ related to the above with data which is
easier to obtain from measurements than the spectral data described.
 \end{abstract}

\section{Introduction and main results}
In this paper we consider some inverse spectral  problems on a closed connected
Riemannian manifold $(M,g)$.  The first motivation to consider
inverse problems on Riemannian  manifolds comes from  spectral
geometry.  The famous problem here, posed by Bochner and formulated by
  Kac in the paper ``Can
one hear the shape of a drum?",  \cite{Kac66}, is the problem of
identifiability of the shape of a $2-$dimensional domain
from the eigenvalues of its Dirichlet Laplacian.
More generally,
 the question is to find the relations between
 the spectrum of a
Riemannian manifold $(M,g)$, i.e. the spectrum of the
Laplace-Beltrami operator $-\Delta_g$ on it,  and  geometry of
this manifold. In particular, one can ask, following Bochner-Kac,
if the spectrum of $-\Delta_g$
 determines the geometry. However, already in 1966, it was known
that,
 in higher dimensions, the answer to this question is negative.
Indeed, in 1964
 Milnor \cite{Mil64} found the first counter-example, a pair
of Riemannian flat tori of dimension $16$ which are isospectral,
i.e. the spectra of their Laplacians coincide, but not isometric.
As for the original Bochner-Kac problem in dimension $2$, the answer
was found only in early 90th. Namely,
in 1985  Sunada \cite{Sun85}
introduced a method of producing examples of non-isometric
isospectral compact connected Riemannian manifolds.
Although in this paper Sunada did not give the answer to the
Bochner-Kac problem, in 1992
Gordon, Webb and Wolpert \cite{GWW92} extended Sunada's method and
settled
in the negative this famous problem by constructing two simply
connected non-isometric plane domains which are isospectral.
Since then there was much activity in this direction with many
beautiful mathematical results regarding relations between the
spectrum and geometry, see e.g.
 \cite{Gor00} and  \cite{Ze} for the current state of art in this area.

It is clear from the above that, in order to determine geometry
of a closed Riemannian manifold, further spectral information is needed.
The nature of this information can be found if we look at
 inverse boundary problems. In this case, the data given is the trace
on $\p M$ of the resolvent of the Laplacian with, say, Neumann
boundary condition. Depending on whether this information is given for
one or many values of the spectral parameter, these inverse boundary
problems were originally posed by Calderon \cite{Cal} and Gel'fand
\cite{Gel}.
These inverse boundary problems were solved, at least
on the level of uniqueness and sometimes conditional stability,
for the Laplace-Beltrami operator and also
the other types of scalar operators,
in  e.g. \cite{AKKL04}--\cite
{BelKur92}, \cite{KatKur98, Kur1, Kur2}, \cite{LasTayUhl03}--\cite{LeeUhl89}, \cite{Nach88}--\cite{
PesUhl05}, \cite{Syl90,SylUhl87}
and monographs \cite{Isa06} or \cite{KKL} with
further references therein.

As the first inverse problem considered in this paper is an analog of the
Gel'fand inverse boundary problem, we first reformulate
the Gel'fand problem in an
equivalent form which, however, has more ``spectral'' flavor.
Namely, let $\lambda_j$ and $\phi_j$ be the eigenvalues and  normalized
eigenfunctions of the Laplace operator with Neumann boundary condition,
 \beq
 \label{Gelfand}
(-\Delta_g-\lambda_j)\phi_j=0\quad\text{in}\quad M,\quad
\p_{\nu}\phi_j|_{\p M}=0;\quad
(\phi_j,\phi_k)_{L^2(M)}=\delta_{jk},
 \eeq
where $\p_{\nu}$ is the normal derivative to $\p M$. Then the Gel'fand
problem
\cite{Gel} is the one of the determination of $(M, g)$ from
the {\it boundary spectral data}, i.e. $\p M,\,
\{\lambda_j, \phi_j|_{\p M}\}_{j=1}^\infty$.
 Note that, due to the Neumann boundary condition in
(\ref{Gelfand}), we do actually know the whole Cauchy data
$\phi_j|_{\p M},\,\p_\nu\phi_j|_{\p M}$ of the eigenfunctions
on $\p M$.
To formulate its
analog for a closed connected manifold, assume as earlier that
we deal with the Laplace operator $-\Delta_g$ on, now closed,
Riemannian manifold $(M,g)$. As earlier, we assume that
our spectral data are given on an $(m-1)-$dimensional
submanifold $\Sigma\subset M, \, \Sigma =\p S$,
where $S \subset M$ is an open set. We note that this type of data
is natural for various physical applications when sources and
receivers are located over some surface in space rather then are
scattered over an $m-$dimensional region or put on, probably remote,
boundary of $M$.
Such localization  is  used e.g. in radars, sonars, and in medical
ultrasound imaging when a single antenna array is used to produce
the wave and to measure the scattered wave. It is typical also
in geosciences/seismology where sources and receivers are often located over
the surface of the Earth or an ocean.

Returning to the possible generalizations of the Gel'fand problem
to inverse problems on closed manifolds and taking into account
that $\Sigma$ splits $M$ into two manifolds with boundary,
$\overline S$ and $M \setminus S$, we believe that the most
straightforward generalization of the Gel'fand problem is
the following:

Assume that we are given  the \emph{Cauchy spectral data},
 \beq
 \label{Cauchy_spectral_data}
\{\Sigma,
(\lambda_j,\phi_j|_{\Sigma},\partial_\nu\phi_j|_{\Sigma})_{j=1}^\infty\}.
\eeq
Does these data determine $(M, g)$?

Then the first main result of the paper is

\begin{thm}
 \label{thm1}
 The Cauchy spectral data (\ref{Cauchy_spectral_data})
determine the manifold $(M, g)$ up to an isometry.
 \end{thm}

Probably a more surprising result
concerns with another
inverse spectral problem associated with $\Sigma$.
In this case we have less data, namely,  the
\emph{Dirichlet spectral data},
\beq
 \label{Dirichlet_spectral_data}
 \{ \Sigma,\
(\lambda_j,\phi_j|_\Sigma)_{j=1}^\infty\},
 \eeq
and would like to determine $(M, g)$ from these data.
It turns out that this is still possible when the set $S$
consists of two subsets,
\[
S= S_1 \cup S_2, \quad {\overline S_1} \cap {\overline S_2} =\emptyset,
\]
if we assume some generic conditions on $S_1, S_2$. To state this result, we define
the Dirichlet spectrum of the Laplace operator on a manifold $\tilde{S}$ with boundary to be a collection of all eigenvalues $\lambda\in \R$ of $-\Delta_g$ such that
\ba
& &-\Delta_g\phi= \lambda\phi,\quad\hbox{in }\tilde S,\quad \phi|_{\p \tilde S}=0\ea
with some non-zero $\phi$.

 \begin{thm}
 \label{thm2} Assume that the Dirichlet spectra of the
Laplace operators on
$\overline{S_1}$, $\overline{S_2}$, $M \setminus S$,
$M \setminus S_1$ and $M \setminus S_2$ are all disjoint. Then
the Dirichlet spectral data determine the manifold up to an
isometry.
 \end{thm}

The paper is organized as follows. In Section \ref{sec_direct} we
present some auxiliary results  for the transmission
problems on $M$ and the metric reconstruction on
$\Sigma$. Section \ref{sec_cauchy} deals with the
reconstruction of the Riemannian manifold $(M,g)$ from the
Cauchy spectral data and contains the proof of Theorem \ref{thm1}.
Section \ref{sec_dirichlet} is devoted to the inverse problem with
the Dirichlet spectral data and provides the proof of Theorem
\ref{thm2}. Section \ref{sec_measurements} contains some further
generalizations of the problem and its alternative formulations
which may be useful in practical applications.

\section{Auxiliary results}

\label{sec_direct}

\subsection{Transmission problem at fixed frequency}
Consider a closed connected smooth $m-$dimensional Riemannian manifold
$(M,g)$ and  an open non-empty set $S \subset M$ with smooth $(m-1)-$dimensional
 boundary $\Sigma:=
 \partial S \neq \emptyset$.
 Let $u:M\to \R$ be a function such that its
restrictions  $u_+,\, u_-$ onto $M_+ = M \setminus S$
and $M_-= \overline S$ are $H^2$-smooth.
 We define the traces of $u$ on the different sides of $\Sigma$ by
\beq
\label{additional}
& &\left(u_{\pm}|_{\Sigma}\right)(x)=\lim_{h\to 0^-}u(x\pm h\nu(x)),\quad x\in\Sigma,\\
\nonumber
& &\left(\partial_{\nu} u_{\pm}|_{\Sigma}\right)(x)=\lim_{h\to 0^-}
\partial_\nu u(x\pm h\nu(x),\quad x\in\Sigma,
 \eeq
where $\nu(x)$ is the unit normal to $\Sigma$ at $x$
pointing towards $M_-$ and the limits are understood in the sense of traces in Sobolev
spaces.
Denote by $[u]=u_+|_{\Sigma}-u_-|_{\Sigma}$
and  $[\p_\nu u] = \p_\nu u_+|_{\Sigma}- \p_\nu u_-|_{\Sigma}$  the jumps
of $u$ and its normal derivative  across
$\Sigma$.

Consider the transmission problem
\begin{equation}
\label{eq_transmission}
\begin{aligned}
&(-\Delta_g-\lambda)u: = - g^{-1/2}\partial_i(g^{1/2}g^{ij}\partial_j u) - \la u
=0 \quad \text{in}\quad M\setminus\Sigma,\\
&[u]=f\quad\text{on}\quad\Sigma,\quad
[\partial_\nu u]=h\quad\text{on}\quad\Sigma, \end{aligned}
\end{equation}
where $g=\det(g_{ij})$ and $[g^{ij}]$ is the inverse matrix of $[g_{ij}]$.
Although, in general, problem (\ref{eq_transmission}) may be considered with
$f \in H^{3/2}(\Sigma), \, h \in H^{1/2}(\Sigma)$, it will be sufficient and convenient
for us to take
 $f,h\in C^\infty(\Sigma)$. Problem (\ref{eq_transmission}) is equivalent to the problem
 \beq
 \label{2.2}
 (-\Delta_g-\lambda)u = h \delta_\Sigma +f \p_\nu \delta_\Sigma,
 \eeq
 where $h \delta_\Sigma$ and $f \p_\nu \delta_\Sigma$ are distributions defined as
 \bfo
 \bra h  \delta_\Sigma, \psi\cet =
 \int_\Sigma h  \psi|_\Sigma \,dS_g, \quad
 \bra f \p_\nu \delta_\Sigma, \psi\cet =
- \int_\Sigma f \left( \p_\nu \psi \right)|_\Sigma \,dS_g,
 \efo
 $dS_g$ being the volume element on $\Sigma$, for any $\psi \in C^\infty(M)$.

 Let $\lambda\not\in \sigma(-\Delta_g)$, where $\sigma(-\Delta_g)$
 is the spectrum of the Laplace operator $-\Delta_g$. Then
 problem (\ref{eq_transmission}) has a unique solution. Moreover, its
 formulation in form (\ref{2.2}) makes it possible to represent this solution,
 $u=u^{f,h}(x,\lambda)$,
 for  $x\in M\setminus\Sigma,$ as a sum of a single- and double-layer potentials,
 \beq
 \label{potentials}
 & &u^{f,h}(x,\lambda) = u^{0,h}(x, \la)+u^{f,0}(x, \la), \\
 \nonumber
 & &u^{0,h}(x, \la)= S_\la h, \,\ u^{f,0}(x, \la)= - D_\la f,  \\
 \nonumber
& &  S_\la h(x)=\int_{\Sigma}
G(x,y;\lambda)h(y)dS_g(y),\quad
D_\la f(x)=\int_{\Sigma}
\partial_{\nu(y)}G(x,y;\lambda)f(y)dS_g(y),
 \eeq
where $G(x,y;\lambda)$ is  Green's
function,
 \[
(-\Delta_g-\lambda)G(x,y;\lambda)=\delta_y(x).
\]
 Note that the
single-layer potential is well-defined on and continuous across   $\Sigma$.
Also,  $D_\la f$
can be
continuously extended from $M_\pm^{int}$ to $M_\pm$ and,
for $x\in \Sigma$,
\beq
\label{double}
 {(D_\la f)_\pm}|_\Sigma(x)=D^{\it o}_\la f(x)\pm
\frac{1}{2}f(x),\quad
D^{\it o}_\la f(x):=\int_{\Sigma}
\partial_{\nu(y)}G(x,y;\lambda)f(y)dS_g(y).
\eeq
 Similar,  the normal derivative of $S_\la h$ is continuously extended, from the left and right,
to $x \in\Sigma$,
\beq
\label{derivative}
{\partial_\nu(S_\la h )_\pm}|_\Sigma (x)= J_\la h(x)
\mp\frac{1}{2}h(x),\quad J_\la h(x):=
\int_{\Sigma}
\partial_{\nu (x)} G(x,y;\lambda) h(y)dS_g(y),
\eeq
 where
the integrals in (\ref{double})-(\ref{derivative}) are weakly singular ones.
Also integrals in (\ref{double})-(\ref{derivative}) are mutually adjoint on $\Sigma$
in the sense of duality
\beq
\label{duality}
& &\bra D^{\it o}_\la f,\,h\cet:= \int_\Sigma \left(D^{\it o}_\la f\right)(x) h(x)\, dS_g(x)
\\
\nonumber &=&
\bra f,\,J_\la h\cet:= \int_\Sigma  f(x) \left( J_\la h\right)(x)\, dS_g(x),
\quad f, h \in C^\infty(\Sigma).
\eeq
(For these results see e.g. \cite{colkres83}. Note that, due to the local nature of constructions in
\cite{colkres83}, they are valid not only for the Euclidian case considered in there
but also for manifolds.)

\subsection{Transmission problem in time domain}
In this subsection we  consider the transmission problem for the wave equation
associated with the Laplace operator, $-\Delta_g$,
 \begin{equation}
 \label{eq_transmission_t}
\begin{aligned}
 &(\p_t^2-\Delta_g)u^{f,h}=0 \quad\hbox{in }(M\setminus \Sigma)\times \R,\\
  &\hbox{[} u^{f,h} \hbox{]} =f \quad\hbox{on }\Sigma\times \R,\quad
  \hbox{[} \partial_\nu u^{f,h}\hbox{]} =h \quad\hbox{on }\Sigma\times \R,\\
  &u^{f,h}|_{t<\min(t_h, t_f)}=0,
\end{aligned}
 \end{equation}
  where $h, f \in C^\infty_+(\Sigma \times \R)$. This space consists of $C^\infty-$ smooth functions
 equal to $0$ for sufficiently large negative $t$, i.e.
 \bfo
 h = 0 \quad \hbox{for}\,\, t<t_h, \quad
f = 0 \quad \hbox{for}\,\, t<t_f.
\efo

The results obtained below will be instrumental,
although for the special case $f=0$, in section \ref{sec_dirichlet}.
However, as considerations for $f=0$ and $f \neq 0$ are parallel, we will consider the general case.

In lemma below, we use the following
spectral decomposition
\beq
\label{Fourier1}
u^{f,h}(x,t)=\sum_{j=1}^\infty u^{f,h}_j(t)\phi_j(x).
\eeq
\begin{lem} (Blagovestchenskii identity) \label{lem_non_stationary}
Given the Cauchy spectral data (\ref{Cauchy_spectral_data})
it is possible to evaluate the Fourier coefficients
$u^{0,h}_j(t)$, $u^{f,0}_j(t)$
of the waves $u^{0,h}(x,t)$, $u^{f,0}(x,t)$, namely,
\beq
\label{time_fourier}
& &
u^{0,h}_j(t)= \int_{t_h}^t \frac{\sin(\sqrt{\la_j}(t-s))}{\sqrt{\la_j}} h_j(s) ds,\quad
h_j(t):= \int_\Sigma h(\cdot,t) \phi_j|_{\Sigma} dS_g, \\ \nonumber
& & u^{f,0}_j(t)=- \int_{t_f}^t \frac{\sin(\sqrt{\la_j}(t-s))}{\sqrt{\la_j}} f_j(s) ds,\quad
f_j(t):= \int_\Sigma f(\cdot,t) \p_\nu \phi_j|_{\Sigma} dS_g,
\eeq
and, for $\la_1=0$, we should substitute $ \frac{\sin(\sqrt{\la_j}(t-s))}{\sqrt{\la_j}}$ by $(t-s)$.
Moreover,
\ba
u^{0,h}(x,t) \in C^\infty_+(\R; H^1(M)), \quad u^{f,0}(x,t) \in C^\infty_+(\R; L^2(M)).
\ea
\end{lem}
\begin{proof}
We  consider first  $u^{0,h}$.
We have
\bfo
u^{0,h}_j(t)=\left(u^{0,h}(x,t), \phi_j(x)\right)_{L^2(M)}=
\int_{M_+}u^{0,h}(x,t) \phi_j(x) dV_g+
\int_{M_-}u^{0,h}(x,t) \phi_j(x) dV_g.
\efo
Thus,   (\ref{eq_transmission_t}) implies that,
\bfo
& &\p_t^2u^{0,h}_j(t)= \int_{M_+}\Delta_g u^{0,h}(x,t) \phi_j(x) dV_g+
\int_{M_-}
\Delta_g u^{0,h}(x,t) \phi_j(x) dV_g= \\
& &
\int_\Sigma \left(\p_\nu u^{0,h}_+-\p_\nu u^{0,h}_-\right)|_\Sigma \phi_j|_\Sigma dS_g -
\int_\Sigma \left(u^{0,h}_+-u^{0,h}_-\right)|_\Sigma \p_\nu\phi_j|_\Sigma dS_g \\
& &+ \int_{M} u^{0,h}(x,t) \Delta_g\phi_j(x) dV_g=h_j(t) -\la_ju^{0,h}_j(t),
\efo
where we use that $\left(\p_\nu u^{0,h}_+-\p_\nu u^{0,h}_-\right)|_\Sigma=h,
\,
\left(u^{0,h}_+-u^{0,h}_-\right)|_\Sigma =0$.
Solving this second order ordinary differential equation together with the initial conditions
$u^{0,h}_j(t_h)=0,\,\p_tu^{0,h}_j(t_h)=0,$ provides the first formula in (\ref{time_fourier}). Similar considerations provide the second
formula in (\ref{time_fourier}).

To prove the second part of Lemma, we rewrite
problem (\ref{eq_transmission_t})
in the form, cf. (\ref{2.2}),
\bfo
\p_t^2 u^{0,h}-\Delta_g u^{0,h} = h \delta_\Sigma,
\quad u^{0,h}|_{t<t_h}=0.
\efo
As, for any $t$, $\,h(\cdot, t) \delta_\Sigma \in H^{-1}(M)$,
we have, for any $k=0,1,\dots,$ that
\beq
\label{H^1_estimate}
\sum_{j=1}^\infty (\la_j+1)^{-1}| \p_t^kh_j(t)|^2 \in C_+^\infty(\R).
\eeq
On the other hand, representation (\ref{time_fourier}) yields that,
for $j\geq 2,$
\bfo
 \p_t^ku^{0,h}_j(t)= \int_{0}^{t-t_h}
 \frac{\sin(\sqrt{\la_j}s)}{\sqrt{\lambda_j}} \p_t^{k}h_j(t-s) ds.
 \efo
 This implies that
\bfo
 \p_t^ku^{0,h}_j(t)= \frac{1}{\la_j}  \p_t^kh_j(t)-\frac{1}{\la_j} \int_{t_h}^t
 \cos(\sqrt{\la_j}(t-s)) \p_s^{k+1}h_j(s) ds,
 \efo
 with obvious modification for $j=1$.
 Thus,
 \bfo
 & &\sum_{j=1}^\infty (\la_j+1)| \p_t^ku^{0,h}_j(t)|^2 \leq
 C\sum_{j=1}^\infty (\la_j+1)^{-1}| \p_t^kh_j(t)|^2 \\
 & &+C'(t-t_h)\int_{t_h}^t
 \sum_{j=1}^\infty (\la_j+1)^{-1}| \p_s^{k+1}h_j(s)|^2\,ds.
 \efo
 This inequality, together with (\ref{H^1_estimate}), implies the desired result for
 $u^{0,h}(x,t)$. Similarly, we prove that
$u^{f,0}(x,t) \in C^\infty_+(\R; L^2(M))$.
\end{proof}
Lemma \ref{lem_non_stationary} immediately implies the following result:
\begin{cor}
\label{trace}
The Dirichlet spectral data (\ref{Dirichlet_spectral_data})
makes it possible to find the trace on $\Sigma \times \R$,
i.e. $u^{0,h}|_{\Sigma \times \R}$, for any
$h \in C^\infty_+(\Sigma \times \R)$.
\end{cor}
\begin{proof}
The result follows from formula (\ref{time_fourier})
taking into the account that, due to
$u^{0,h} \in C^\infty_+(\R; \, H^1(M))$,
the series (\ref{Fourier1}) converges, for any $t$,
in $H^1(M)$. Therefore, the trace
\bfo
u^{0,h}|_{\Sigma \times \R}(x,t)=\sum_{j=1}^\infty u^{0,h}_j(t)
\phi_j|_{\Sigma}(x),
\efo
where the right-hand side converges in
$C^\infty_+(\R; \, H^{1/2}(\Sigma))$
\end{proof}

Let us now warm up by considering a simple inverse problem when the part
$M_-$ of the manifold $M$  is  known and show
how Lemma \ref{lem_non_stationary} and Corollary \ref{trace}
can be used to recover $M_+$.

 \begin{exmp}  Assume
that, in addition to (\ref{Dirichlet_spectral_data}),
we know the manifold $(\overline{M}_-, g|_{\overline{M}_-})$.
 Then the manifold $(M_+,g|_{M_+})$ can be recovered up to an isometry.

 Indeed, let
 $h \in C^\infty_+(\Sigma \times \R_+)$, i.e.
 $h \in C^\infty_+(\Sigma \times \R)$ and $t_h >0$.
Using Corollary
 \ref{trace}, we can find
 $u^{0,h}|_{\Sigma\times\R_+}$.

 Consider now the Dirichlet initial-boundary value
 problem in $M_- \times \R_+$ with Dirichlet data being
 $u^{0,h}|_{\Sigma\times\R_+}$. As $(M_-, \,g|_{M_-})$
 is known, we can then find $u^{0,h}|_{M_- \times\R_+}$
 and, therefore,
  the normal derivative $\p_\nu u^{0,h}_-|_{\Sigma\times\R}$. Thus, we can find $\p_\nu u^{0,h}_+|_{\Sigma\times\R}=\p_\nu
 u^{0,h}_-|_{\Sigma\times\R}+h$.

It is shown in the proof of Lemma \ref{Energy} below
that, when $h$ runs over $C^\infty_+(\Sigma \times \R_+)$,
then $u^{0,h}|_{\Sigma\times\R_+}$ also runs over the whole
$C^\infty_+(\Sigma \times \R_+)$. Therefore, the set of pairs,
  \[
 \{(u^{0,h}|_{\Sigma\times\R},\p_\nu u^{0,h}_+|_{\Sigma\times\R})\,:\,  h\in C^\infty_+(\Sigma \times \R)\},
 \]
defines the graph of the non-stationary
 Dirichlet-to-Neumann map for $M_+$.

 It is, however, shown in e.g. \cite{KKL}
that this map determines  $(M_+,g|_{M_+})$ upto an isometry.
 \end{exmp}

%Next we return to the case when both parts $M_+$ and $M_-$ of a
%manifold $M$ are unknown.

\subsection{Metric on $\Sigma$.}
We complete this section showing that
the Dirichlet spectral data determines the metric on $\Sigma$.
This result will be needed later, in sections \ref{sec_cauchy}
and \ref{sec_dirichlet}.

\begin{lem} \label{lem_volume}
 The Dirichlet spectral data (\ref{Dirichlet_spectral_data})
 determine the distance
function on $\Sigma$ and, therefore, the inherited metric tensor,
$g|_{\Sigma}$,  on
$\Sigma$.
 \end{lem}

\begin{proof} Let $H(x,y;t)$ be the  heat kernel,
 \begin{align*}
 \p_tH(x,y;t)-\Delta_gH(x,y;t)=\delta_y(x)\delta(t),\quad x,y\in M; \quad
 H|_{t=0}=0.
 \end{align*}
Then, uniformly  in $M\times M$,
\bfo
 t\log H(x,y;t)\to -\frac{1}{4}d_M(x,y),
\efo
 as $t\to 0$, see e.g.
\cite{var}.
On the other hand, for $t>0$,
\bfo
 H(x,y;t) = \sum_{j=1}^\infty e^{-\la_j t} \phi_j(x) \phi_j(y),
 \efo
 where the convergence takes place in $C^\infty(M \times M \times (0, \infty))$.
Thus, we can determine the distance $d_M(x,y)$
between any points $x,y\in\Sigma$ using the Dirichlet
spectral data (\ref{Dirichlet_spectral_data}).
Then the distance along $\Sigma$ is
given by
 \[
 d_\Sigma(x,y)=\lim_{\varepsilon\to 0}\inf\sum_{j=0}^{N-1}d_M(z_j,z_{j+1}),
 \]
where the infimum is taken over all finite sequences of points
$z_0,z_1,\dots,z_N\in\Sigma$ such that $z_0=x$, $z_N=y$ and
$d_M(z_j,z_{j+1})\le \varepsilon$. Having at hand
$d_\Sigma(x,y)$, we
can determine the inherited metric tensor, $g|_\Sigma$ on $\Sigma$.
\end{proof}

\section{ Inverse  problem with Cauchy spectral data. Theorem \ref{thm1}.}

\label{sec_cauchy}
In this section we  develop a procedure to reconstruct the
Riemannian manifold $(M,g)$ from the \emph{Cauchy spectral data}
(\ref{Cauchy_spectral_data}).

\subsection{From Cauchy spectral data to the response operator}
Let us consider the transmission problem
(\ref{eq_transmission}) with $\lambda\not\in \sigma(-\Delta_g)$.
Denoting, as in section \ref{sec_direct}, its solution by $u^{f,h}(x, \la)$,
we define the \emph{response operator} by
setting
 \ba R_\lambda(f,h)=u^{f,h}_+(\la)|_{\Sigma}.
 \ea

\begin{lem} \label{lem_potentials}
 The Cauchy spectral data (\ref{Cauchy_spectral_data})
 determine the response operator $R_\lambda$ for all $\lambda\notin\sigma(-\Delta_g)$.
 \end{lem}

\begin{proof}
Note that Green's function has the following
spectral representation
 \beq
 \label{3}
G(x,y;\lambda)=\sum_{j=1}^{\infty}\frac{\phi_j(x)\phi_j(y)}{\lambda-\lambda_j},
 \eeq
 where the sum converges the  sense of operators in $L^2(M)$ and we
 assume, without loss of generality,
 that the eigenfunctions $\phi_j$ are real-valued.
 Thus,  given the Cauchy spectral data and using (\ref{3}), one can
 formally construct the single- and double-layer potentials on
 $\Sigma$,
 \begin{equation}
 \label{eq_sigle_green}
 S_\la h(x)|_{\Sigma}\underset {\text{formally}} =\int_{\Sigma}
 \big[\sum_{j=1}^\infty\frac{\phi_j(x)|_{\Sigma}\phi_j(y)}{\lambda-\lambda_j}\big]h(y)dS_g(y),
 \end{equation}
 \begin{equation}
 \label{eq_double_green}
D_{\la}^{\it o} f(x)\underset {\text{formally}} =\int_{\Sigma}
 \big[\sum_{j=1}^\infty\frac{\phi_j(x)|_{\Sigma}\partial_{\nu(y)}\phi_j(y)}{\lambda-\lambda_j}\big]f(y)dS_g(y).
 \end{equation}
 Since
 \ba
 R_\lambda(f,h)=S_\la h-(\frac 12+D_{\la}^{\it o})f,
 \ea
 it looks that data (\ref{Cauchy_spectral_data}) directly determines $R_\la$.
However, we face the difficulty that  series
 (\ref{3}) does not converge pointwise. To deal with this difficulty, consider first the case when
 $f=0$. To determine  the coefficients in the Fourier expansion
 \begin{equation}
 \label{eq_u}
 u^{0,h}(x,\,\lambda)=\sum_{j=1}^\infty (u^{0,h}(\lambda),\phi_j)_{L^2(M)}\phi_j(x),
 \end{equation}
 %converges in $H^1(M)$ as well as in $L^2(M)$.
we use Green's formula
to get
   \beq
   \label{eq_sol_u}
& &  (u^{0,h},\phi_j)_{L^2(M)}= -\frac{1}{\lambda-\lambda_j}
(\int_{M_+}+\int_{M_-}) (\Delta_g
u^{0,h}\phi_j-u^{0,h}\Delta_g\phi_j)dV_g=
\\ \nonumber
& &
-\frac{1}{\lambda-\lambda_j}
\int_{\Sigma}(\partial_\nu u^{0,h}_+\phi_j|_{\Sigma}-\partial_\nu
u^{0,h}_-\phi_j|_{\Sigma}-
u^{0,h}_+\partial_\nu\phi_j|_{\Sigma}+u^{0,h}_-\partial_\nu\phi_j|_{\Sigma})dS_g=
-\frac{1}{\lambda-\lambda_j}
\int_{\Sigma}h\phi_j|_{\Sigma}dS_g.
 \eeq
Since the series (\ref{eq_u}), (\ref{eq_sol_u}) converges
% \begin{equation}
% \label{eq_sol_u}
% u^{0,h}(x,\lambda)=-\sum_{j=1}^\infty
%\big[ \int_{\Sigma} \phi_j(y)
% h(y)dS_g(y) \big]  \frac{\phi_j(x)}{\lambda-\lambda_j},
%  \end{equation}
 in $H^1(M)$, so that
 the trace is given by
 \[
 u^{0,h}(x,\lambda)|_{\Sigma}=-
 \sum_{j=1}^\infty
\big[ \int_{\Sigma}  \phi_j(y) h(y)dS_g(y) \big]
\frac{\phi_j(x)|_{\Sigma}}{\lambda-\lambda_j},
  \]
   where the series
converges in $H^{1/2}(\Sigma)$. Hence \eqref{eq_sigle_green} is
well-defined. To compute \eqref{eq_sigle_green} we  also need to know
 the Riemannian volume $dS_g(x)$ of $\Sigma$. By
 Lemma
\ref{lem_volume}, it can be found from data (\ref{Cauchy_spectral_data}).

Let us now show that \eqref{eq_double_green} is well-defined.
First note that $u^{0,h}|_{M_+}\in H^1(M_+) \cap H^2_{loc}(M_+)$. Thus, we can define,
for any $\e >0$,
the normal derivative $\p_\nu u^{0,h}|_{\Sigma^\e}$,
where $\Sigma^\e = \{x \in M_+: d(x, \Sigma)=\e\}$. Clearly, for any
$\Psi\in H^{1}(M_+)$,
 \beq
 \label{Q} \nonumber
& &  \int_{\Sigma^\e} (\p_\nu
u^{0,h}|_{\Sigma^\e})\Psi|_{\Sigma^\e}\, dS_g= \int_{M_+^\e} \Delta_g
u^{0,h} \,\Psi\, dV_g + \int_{M_+^\e} (\nabla_g
u^{0,h},\, \nabla_g \Psi)_g\,dV_g
\\
& &
=- \la \int_{M_+^\e}
u^{0,h} \,\Psi\, dV_g + \int_{M_+^\e} (\nabla_g
u^{0,h},\, \nabla_g \Psi)_g\,dV_g
 \eeq
where $M_+^\e= \{x \in M_+: d(x, \Sigma) \geq\e\}$
and we have used equation (\ref{eq_transmission}).
As the right-hand side of (\ref{Q}) has a limit, when $\e \to 0$, and $\Psi|_{\Sigma}$
runs over $H^{1/2}(\Sigma)$ when $\Psi$ runs over $H^1(M_+)$,
this  defines $\p_\nu
u^{0,h}_+|_{\Sigma}$. As we can choose, for any $\Psi|_{\Sigma}
\in H^{1/2}(\Sigma)$ its extension $\Psi$ so that
$\|\Psi\|_{H^1}\leq C\|\Psi|_{\Sigma}
\|_{H^{1/2}}$,
 \beq
 \label{eq_u_trace}
 \|\p_\nu
u^{0,h}(\lambda)_+|_{\Sigma}\|_{H^{-1/2}(\Sigma)}\leq
C\|u^{0,h}(\lambda)|_{M_+}\|_{H^1(M_+)}.
 \eeq

Define the sources-to-Dirichlet operator by setting
 \ba
& &J_\lambda h=\p_{\nu}u^{0,h}(\lambda)_+|_{\Sigma},
 \ea
 see (\ref{derivative}).
Taking $\lambda$-derivative of \eqref{eq_sol_u}, we get
 \ba
 \frac{\p}{\p \lambda} u^{0,h}(x, \lambda)=
 \sum_{j=1}^\infty  \big[ \int_\Sigma \phi_j(y) h(y)\,dS_g(y)
 \big] \frac{\phi_j(x)}{
(\lambda_j-\lambda)^{2}}
=
(-\Delta_g-\lambda)^{-1}u^{0,h}(\lambda)
 \ea
 that converges in
$H^{3}(M).$ So we have a well-defined object
 \ba
 \frac{\p}{\p \lambda} \p_\nu u^{0,h}(\lambda)|_{\Sigma}
=\sum_{j=1}^\infty \big[
\int_\Sigma \phi_j(y)  h(y)\,dS_g(y) \big]
\frac{\p_{\nu(x)}\phi_j(x)|_\Sigma
}{ (\lambda_j-\lambda)^{2}} ,
\label{4B}
 \ea where the convergence holds in
$H^{3/2}(\Sigma)$.
As
$
 \frac{\p}{\p \lambda} \left(J_\lambda h \right)=
\frac{\p}{\p \lambda} \p_\nu u^{0,h}(\lambda)|_{\Sigma},
 $
we can compute,  for any $h\in C^\infty(\Sigma)$,
$\frac{\p}{\p \lambda} \left(J_\lambda h \right)$ using
 the Cauchy spectral data of
$\Sigma$.

Let $\lambda\in \R$, $\lambda\not=\lambda_j$, and let
$\gamma_T\subset \C$ be the line segment from $\lambda$ to $iT$.
As
 \ba J_\lambda h=\int^\lambda_{iT} \frac{\p}{\p
\tau} (J_\tau h)\,d\tau+ J_{iT} h,
 \ea
 we have
 \bfo
 J_\lambda
h=\lim_{T\to\infty}\big(\int^\lambda _{iT} \frac{\p}{\p \tau}
(J_\tau h)\,d\tau+ J_{iT} h\big).
 \efo

By Lemma \ref{lem_asymptotics} below, we get $
\lim_{T\to\infty} J_{iT} h= \lim_{T\to\infty} \p_\nu u^{0,h}(iT)_+|_{\Sigma}=\frac12 h$. This implies that
\ba
J_\lambda
h=\frac12 h + \lim_{T\to\infty}\int^\lambda _{iT} \frac{\p}{\p \tau}
(J_\tau h)\,d\tau,
\ea
where the right-hand side  can be computed using
the Cauchy spectral data.

To complete the proof, we recall,
see equation (\ref{duality}),
that $J_\lambda$ is adjoint  of  $D_\la^{\it 0}$. Thus
 we can find
$D_\la^{\it o} f$ using the Cauchy spectral data.
\end{proof}

In the proof of the above Lemma we used the following asymptotics, with respect to singularity,
of $u^{0,h}$ near $\Sigma$.

\begin{lem}
\label{lem_asymptotics}
Let $\lambda\in\mathbb{C}$, $\Im\lambda>0$ and $\emph{\Re}{\sqrt{-\la}} < 0$. Then
uniformly for
$|\arg(\la)| \geq \delta$,
\beq
\label{traces}
u^{0,h}(\la)_\pm|_\Sigma \rightarrow 0,\quad
\p_\nu u^{0,h}(\la)_\pm|_\Sigma \rightarrow \pm \frac12 h,
\quad \hbox{as} \,\, \la \to \infty,
\eeq
in $H^{3/2}(\Sigma),\,H^{1/2}(\Sigma)$, correspondingly.
\end{lem}

\begin{proof} Let us first fix local coordinates on $M$ near $\Sigma$, $x=(x',x^m)$, where
 $x'=(x^1,\dots,x^{m-1})$ are some local coordinates  on $\Sigma$ and $x^m$ is the signed distance to $\Sigma$,
\bfo
x^m= \pm \hbox{dist}(x, \Sigma) \quad \hbox{for} \,\, x \in M_\pm.
\efo
In these coordinates, we introduce
\begin{equation}
\label{singularity}
\begin{aligned}
v^{0,h}_\pm(x; \la):= \begin{cases} \frac{h(x')}{2\sqrt{-\la}}  e^{\,x^m {\sqrt{-\la}} } \zeta(x^m), & x\in M_+,\\
\frac{h(x')}{2\sqrt{-\la}}  e^{-x^m {\sqrt{-\la}} } \zeta(x^m), & x\in M_-,
\end{cases}
\end{aligned}
\end{equation}
where
$\zeta(x^m)$ is a smooth cut-off function equal to $1$ near $x^m=0$ supported in $(-a,a)$ with sufficiently small $a>0$. Outside the a-neighborhood of $\Sigma$, the functions $v_\pm$ are defined to be zero.  Writing
\[
\Delta_g=\p^2_{x^m}+p(x)\p_{x^m}+Q(x,\p_{x^1},\dots,\p_{x^{m-1}})
\]
in the above coordinates and using
the fact that
\[
\int_0^a|e^{\,x^m\sqrt{-\lambda}}\zeta(x^m)|^2dx^m\le C(\sqrt{-\lambda})^{-1},
\]
we see that (\ref{2.2}), (\ref{singularity}) yield
\bfo
u^{0,h}(x, \la) =v^{0,h}(x, \la)+w^{0,h}(x, \la),
\efo
where $w^{0,h}\in H^2(M)$ satisfies
\bfo
& &(-\Delta_g-\lambda) w^{0,h}= {\it H}^h(\la),
\quad
 \|{\it H}^h(\la)\|_{L^2(M)} \leq C_h(1+|\la|)^{-1/4}.
\efo

As
$\|(-\Delta_g-\lambda)^{-1}\| \leq \hbox{dist}(\la, \sigma(-\Delta_g))$, where the norm is
the operator norm in $L^2(M)$, this implies
\bfo
\|w^{0,h}(\la)\|_{H^2(M)} \leq C_{h,\delta} |\la|^{-1/4},
\efo
when $|\hbox{arg}(\la)| \geq \delta >0$ and $|\lambda|>1$.
Combining this estimate with (\ref{singularity}), we see \eqref{traces}.

\end{proof}

\begin{remark}
Analyzing the behaviour of $u^{f,h}$ near $\Sigma \times \R$, we can show, cf. considerations
leading to (\ref{traces}) and (\ref{singularity}), that
\beq
\label{smoothness}
u^{f,h}(x,t)|_{M_\pm \times \R} \in C^\infty_+(M_\pm \times \R),
\eeq
meaning that $u^{f,h}(x,t)|_{M_\pm^{int} \times \R}$ may be continued to
$M_\pm \times \R$ to satisfy (\ref{smoothness}).
\end{remark}

\subsection{Reconstruction of the manifold using the response operator}
Recall, see e.g. \cite{KKL}, that, if $(N,g),\, \p N \neq \emptyset$, then its
Neumann-to-Dirichlet operators, $\Lambda_\la(N)$ are defined as
\bfo
\Lambda_\la(N) \psi= w^{\psi}(\la)|_{\p N},
\efo
where $w^{\psi}(x, \la)$ is the solution to the Neumann problem
\bfo
-\Delta_g w^{\psi}(x, \la) =\la w^{\psi}(x, \la), \,\, x \in N^{int}, \quad
\p_\nu w^{\psi}(x, \la)|_{\p N}=\psi,
\efo
for $\la \notin \sigma(-\Delta_g^N), \, \sigma(-\Delta_g^N)$ being the spectrum
of the Neumann Laplacian on $N$.
\begin{lem}
 \label{lem_neumann-to-dirichlet}
 Given the  Cauchy spectral data (\ref{Cauchy_spectral_data})
  it is possible to find the Neumann-to-Dirichlet
operators $\Lambda_\la(M_\pm)$ for $\la \notin \sigma(-\Delta_\pm^{N})$,
where $-\Delta_\pm^{N}$ stands for the Neumann Laplacian on $M_\pm$.
 \end{lem}
\begin{proof}
We start  with $\Lambda_\la(M_-)$, assuming $\la \notin
\left(\sigma(-\Delta_g) \cup\sigma(-\Delta^N_-) \cup
\sigma(-\Delta^D_+)\right)$, where $-\Delta^D_\pm$ is the Dirichlet Laplacian in $M_\pm$,
correspondingly. Then, for any  $h \in C^\infty(\Sigma)$,
there is a unique solution, $w^h_-(x,\la) \in C^\infty(M_-)$, satisfying
\bfo
-\Delta_g w^h_-(x,\la) =\la w^h_-(x,\la) \,\, \hbox{in} \,\, M_-,\quad
\p_\nu w^h_-(\cdot,\la)|_{\Sigma}= -h,
\efo
where, as in equation (\ref{additional}),  $\nu$ is the unit normal pointing towards $M_-$.

Consider
\[
w^h(x, \la)=\left\{ \begin{array}{ccc}
0 &  \text{in} & M_+,\\
w^h_-(x,\la) &\text{in} & M_-.
     \end{array}\right.
\]
 Clearly, $w^h(x, \la)$ solves (\ref{eq_transmission}) with
 $[w^h]:=f=-w^h_-(\cdot,\la)|_{\Sigma},\, [\p_\nu w^h]=h$.  Moreover, with this $f$ and $h$,
 \beq
 \label{defining_equation}
 R_\la(f,h)=0.
 \eeq
 These considerations show that, for any $h$, there is $f$ such that
(\ref{defining_equation})
 is satisfied and we can consider (\ref{defining_equation}) as an equation
for $f$  when $h$ is given.
 Let us show that the solution to (\ref{defining_equation}) is unique if
 $\la \notin \left(\sigma(-\Delta_-^N) \cup \sigma(-\Delta_+^D) \right)$.
 This will allow us to uniquely
 define
 $f=f^h(\la)=-w^h_-(\cdot,\la)|_{\Sigma}$ as the solution
to (\ref{defining_equation}).
Then,
 \bfo
 \Lambda_\la(M_-) h= - f^h(\la).
 \efo
 To prove uniqueness, assume that there is $f$ such that
\bfo
R_\la(f,0)=0.
\efo
As  $\la \notin  \sigma(-\Delta_+^D)$, this implies that $u_+^{f,0}(x, \la)=0$.
As $[\p_\nu u^{f,0}]=0$, we see that
\bfo
\p_\nu u_-^{f,0}(\cdot, \la)|_{\Sigma}=0.
\efo
However, $\,\la \notin \sigma(-\Delta_-^N) $, so that $ u_-^{f,0}(x, \la)=0$, i.e. $f=0$.

Combining with Lemma \ref{lem_potentials},
 we see that the Cauchy spectral data (\ref{Cauchy_spectral_data})
  determine $\Lambda_\la(M_-)$ for
$\la \notin
\left(\sigma(-\Delta_g) \cup \sigma(-\Delta^N_-) \cup
\sigma(-\Delta^D_+)\right)$. Since $\Lambda_\la(M_-)$ is a meromorphic
operator-valued function with simple poles at $\sigma(-\Delta^N_-)$, this determines
$\Lambda_\la(M_-)$ uniquely.

As $u^{f,h}_-|_{\Sigma}= R_\la(f,h) -f$, we can repeat
the previous arguments for $\Lambda_\la(M_+)$.
\end{proof}

Theorem \ref{thm1} follows from Lemma \ref{lem_neumann-to-dirichlet} taking into account that
$\Lambda_\la(M_\pm)$ determine $(M_\pm, g_\pm)$ up to an isometry, see
\cite{KKL}, section 4.1. Thus to recover $(M,g)$ we should just glue
$(M_-,  g_-)$ and $(M_+,  g_+)$ along given $\Sigma$.

\section{Inverse problem with Dirichlet spectral data. Theorem \ref{thm2}.}

\label{sec_dirichlet}
In this section, we will develop a procedure to reconstruct the
Riemannian manifold $(M,g)$ from the \emph{Dirichlet spectral data}
(\ref{Dirichlet_spectral_data}).
We will assume that $S \subset M$
consists of two open subsets $S_1,S_2, \, \overline{S_1} \cap \overline{S_1} =\emptyset,
\, S= S_1 \cup S_2$. As in section \ref{sec_cauchy}, we assume that $\Sigma:= \p S=\Sigma_1 \cup \Sigma_2,\,\Sigma_i=
 \p S_i, i=1,2,$ are smooth. Moreover, we assume that the spectra
$\sigma(-\Delta^D( \overline{S_i})),\,  \sigma(-\Delta^D( M \setminus S_i)), i=1,2,$
and $\sigma(-\Delta^D(M \setminus S))$ are all disjoint.

\subsection{An approximate controllability result.}
Consider the following transmission problem
 \begin{equation}
  \label{eq_problem_f=0}
\begin{aligned}
 &(\p_t^2-\Delta_g)u=0,\quad\hbox{in }(M\setminus \Sigma)\times \R, \\
 &\hbox{[} u \hbox{]}_{\Sigma} =0,\quad
 \hbox{[} \p_\nu u\hbox{]}_{\Sigma} =h\in C^{\infty}_+(\Sigma\times \R),\\
 &u|_{t<t_h}=0,
\end{aligned} \end{equation}
 and denote by $u(x,t)=u^{0,h}(x,t)$ its
solution. Note that problem (\ref{eq_problem_f=0}) coincides with problem
(\ref{eq_transmission_t}) with $f=0$.

By Lemma \ref{lem_non_stationary},
 $ u^{0,h}\in C^\infty_+(\R; H^{1}(M))$ and we can define an operator
\[
W: C^{\infty}_+(\Sigma\times \R)\to H^1(M),\quad
W h:=u^{0,h}(0),
\]
which is called the wave operator associated with problem (\ref{eq_problem_f=0}).

\begin{thm} \label{lem_controll}
Let $\sigma(-\Delta^D(\overline S)) \cap \sigma(-\Delta^D(M \setminus S))
= \emptyset$. Then the
 set
 \beq
\label{Y}
Y=\{W h:\, h\in
C^{\infty}_+(\Sigma\times \R)\}
 \eeq
  is dense in $H^1(M)$.
\end{thm}

\begin{proof} Assume that $\psi\in (H^1(M))'=H^{-1}(M)$ is orthogonal to $Y$,
 \beq
 \label{orthogonal}
 (u^{0,h}(0),\psi)_{H^1(M)\times H^{-1}(M)}=0
 \eeq
for all $h\in C^\infty_+(\Sigma\times \R)$.
  Let $e$ be the solution to the problem,
\beq\label{2.1}
& &e_{tt}-\Delta_g e = 0,\quad \hbox{ in }M\times \R,\\
& &e|_{t=0} =0,\  e_t|_{t=0}  = \psi. \nonumber
\eeq
Then similar considerations to those at the end of proof
of Lemma \ref{lem_non_stationary} show that  the weak solution,
$e(x,t)$ of \eqref{2.1} satisfies,
 \begin{equation}
  \label{eq: less regular energy estimate A}
  e(x,t)=\sum_{j=2}^\infty \frac{\sin(\sqrt {\lambda_j}t)}
{\sqrt {\lambda_j}}(\psi,\phi_j)_{H^{-1}(M)\times H^1(M)}\phi_j+t(\psi,\phi_1)_{H^{-1}(M)\times H^1(M)}\phi_1
 \end{equation}
 and $e\in
C(\R;
L^2(M))\cap C^1(\R; H^{-1}(M))$. Observe that, as $\|\phi_j|_\Sigma \|_{H^{1/2}} \leq C(\la_j+1)^{1/2}$
 and $\la_j > C j^{2/m}$, $e$ has a well-defined trace, in ${\mathcal S}'(\Sigma\times\R)$,
 on $\Sigma\times\R$ with
 \beq
 \label{trace_sum}
\sum_{j=2}^J \frac{\sin(\sqrt {\lambda_j}t)}
{\sqrt {\lambda_j}}(\psi,\phi_j)\phi_j|_{\Sigma}+ t(\psi,\phi_1)\phi_1|_{\Sigma} \rightarrow e(x,t)|_{\Sigma\times\R}, \quad
\hbox{as}\,\, J\to \infty,
\eeq
in ${\mathcal S}'(\Sigma\times\R)$.

Let us show that $e(x,t)|_{\Sigma\times\R}=0$. Choosing $h \in C^\infty_0(\Sigma \times \R_-)$
and using
 Green's formula, we obtain from (\ref{eq_problem_f=0})--(\ref{2.1}) that
 \ba 0&=&\int_{M\times
\R_-}[u^{0,h} {(e_{tt}-\Delta_g e)} -(u^{0,h}_{tt}-
\Delta_g u^{0,h})  {e}]\,dV_g\,dt\\
 &=&(u^{0,h}(T),  {\psi})_{H^1(M)\times H^{-1}(M)}\,+
  \int_{\Sigma \times\R_-} h\, {e}\,dS_g\,dt\,=
 \int_{\Sigma \times\R_-} h\, { e}\,dS_g\,dt.
 \ea
This yields that $
e|_{\Sigma \times (-\infty,0)}=0. $ As by (\ref{2.1})
 $e(x,s)=-e(x,-s)$, we see that
 \beq
 \label{eq: C data}
 \hbox{supp}\left(e|_{\Sigma
\times \R}\right)=\Sigma \times\{0\}.
 \eeq

 Next we show that
 \beq
 \label{trace_of_e}
 e|_{\Sigma \times (-1, 1)} \in \tilde H^{-1/2}(\Sigma \times (-1, 1))
 :=\left(H^{1/2}(\Sigma \times (-1, 1))\right)'.
 \eeq

 Let $X$ be a local, near $\Sigma$, vector field on $\overline S$ such that
 $
 X|_\Sigma =\p_\nu|_\Sigma.
 $
 Let $h\in H^{1/2}(\Sigma\times (-1,1))$ and $H\in H^1(S\times (-1,1))$ be its continuation into $S \times(-1, 1)$, such that $H=0$ outside the domain of definition of $X$ and
 \bfo
  \|H\|_{H^1}\leq C\|h\|_{H^{1/2}}.
  \efo
  Denote by $E$ the primitive, with respect to $t$, of $e$ in $S \times (-1, 1)$,
 \begin{equation}
  \label{E}
 E(x,t) = -\sum_{j=2}^\infty \frac{\cos(\sqrt {\lambda_j}t)}
{\lambda_j}(\psi,\phi_j)\phi_j(x)+ \frac{t^2}{2}(\psi,\phi_1)\phi_1(x),
\end{equation}
$ E \in C(\R; H^{1}(S)) \cap C^1(\R; L^{2}(S))$. Integrating by
parts, we get
 \beq
 \label{embedding}
& & \int_{\Sigma\times (-1,1)} h\, e\,dS_g\,dt
=
 \int_{S\times (-1,1)} (H \cdot Xe- X^cH\cdot e)\,dV_g\,dt\\ \nonumber
& & = \int_S \left([H \cdot XE]|_{t=1}- [H \cdot XE]|_{t=-1} \right) \,dV_g -
  \int_{S\times (-1,1)} (\p_tH \cdot XE+ X^cH\cdot e)\,dV_g\,dt,
  \eeq
where $X^c$ is the first-order operator adjoint
to $X$.
By (\ref{eq: less regular energy estimate A}), (\ref{E}), the right-hand side of
(\ref{embedding})  can be estimated
by
\bfo
C\|H\|_{H^1(S \times (-1, 1))} \leq C' \|h\|_{H^{1/2}(\Sigma \times (-1, 1))}.
\efo
Thus the left-hand side of (\ref{embedding}) is bounded for any $h \in  H^{1/2}(\Sigma \times (-1, 1))$,
proving (\ref{trace_of_e}).

 Now \eqref{eq: C data} implies that
$
e(x,t)|_{\Sigma
\times \R}=\sum_{i=0}^I e_i(x) \partial^i\delta(t),
$
with some finite $I$,
 see e.g. \cite[ex. 5.1.2]{Hor03}. Thus, \eqref{trace_of_e} yields
 \beq
 e(x,t)|_{\Sigma \times \R}=0.
 \label{zero}
 \eeq

 The last step of the proof is to show that this equation yields that $e=0$ in $M \times \R$.
Using relation (\ref{trace_sum}),
equation (\ref{zero}) and making the partial Fourier transform, $t \to k$, we see that
the distribution, $\hat e(x,k) \in {\mathcal S}'(\Sigma \times \R)$, satisfies
\ba
\hat e(x,k) &=& i\big(\frac12 \sum_{j=2}^\infty \frac{\delta(k-\sqrt{\la_j})}{\sqrt{\lambda_j}}
(\psi,\phi_j)\phi_j|_{\Sigma}\\
&&
+\delta'(k) (\psi,\phi_1)\phi_1|_{\Sigma}
 -
\frac12 \sum_{j=2}^\infty \frac{\delta(k+\sqrt{\la_j})}{\sqrt{\lambda_j}}
(\psi,\phi_j)\phi_j|_{\Sigma}\big)
=0.
\ea
This implies that $(\psi,\phi_1)_{H^{-1}(M)\times H^1(M)}=0$ and,
for any $\tilde j =2,  \dots$,
\beq
\label{dirichlet_S}
 \sum_{\la_j=\la_{\tilde j}}  (\psi,\phi_j)_{H^{-1}(M)\times H^1(M)}\phi_j|_{\Sigma}=0,
\eeq
where the last sum takes into account eigenspaces of
an arbitrary multiplicity.
Consider the function
\bfo
\Phi(x):= \sum_{\la_j=\la_{\tilde j}}  (\psi,\phi_j)_{H^{-1}(M)\times H^1(M)}\phi_j(x),
\quad x \in M.
\efo
It satisfies the Dirichlet boundary condition,
$\Phi|_{\Sigma}=0$, see (\ref{dirichlet_S}), and, as $\phi_j$ are
eigenfunctions of $-\Delta_g$ with $\la_j=\la_{\tilde j}$, the equation
\bfo
-\Delta_g \Phi(x)= \la_{\tilde j} \Phi(x), \quad x \in M.
\efo
Thus $\Phi|_S$ is an eigenfunction of $-\Delta^D(S)$,
while $ \Phi|_{M\setminus S}$
is an eigenfunction of $-\Delta^D(M\setminus S)$.
However, as
$\sigma(-\Delta^D(\overline S)) \cap \sigma(-\Delta^D(M \setminus S))
= \emptyset$, we have that $ \Phi|_{\overline S}=0$ or $ \Phi|_{M \setminus S}=0$. In any case,
by the  uniqueness of zero-continuation for elliptic equations, this yields that
$\Phi=0$ everywhere in $M$.  As different $\phi_j$, corresponding
to $\la_j=\la_{\tilde j}$, are linearly independent, this implies that
$(\psi,\phi_j)_{H^{-1}(M)\times H^1(M)}=0$ for all $j=1,2, \dots$. Thus,
$e=0$ in $M \times \R$ and, therefore, $\psi=0$.
\end{proof}

\subsection{Approximate controllability with given trace at final
time}
In this section we denote $\tilde \Sigma$ to be either $\Sigma_i, i=1,2,$
or$\,\Sigma$.
Lemma \ref{lem_non_stationary} makes it possible to introduce a quasinorm
\beq
\label{quasinorm}
|h|^2:= \|Wh \|^2_{H^1(M)} =
\sum_{j=1}^\infty(\la_j+1)|u_j^{0,h}(0)|^2.
\eeq
It is classical for the control theory, see e.g.
\cite{KL-Dirac} or \cite{LLT} in the context of inverse problems,
to define the space ${\it D}^1$ of the generalized sources
by introducing the equivalence relation,
\bfo
h \equiv_{\it E} \tilde h \quad \hbox{if} \,\, u^{0,h}(0)=
u^{0,{\tilde h}}(0),
\efo
and completing $C^\infty_+(\Sigma  \times \R)/{\it E}$ with respect to (\ref{quasinorm}),
\bfo
{\it D}^1:= \hbox{cl}\left(C^\infty_+(\Sigma  \times \R)/{\it E} \right).
\efo
Then, by Theorem \ref{lem_controll}, we can extend the wave operator
$W$, see (\ref{Y}),
from $C^\infty_+(\Sigma  \times \R)$ onto ${\it D}^1$,
\bfo
Wh:=u^{0,h}(0), \quad W:{\it D}^1 \rightarrow H^1(M),
\efo
as a unitary operator.

Moreover, as $Wh=\sum_{j=1}^\infty \kappa_j^h \phi_j \in H^1(M)$ and the Fourier coefficients
$\kappa_j^h$,
for any $h \in {\it D}^1, $ can be explicitly evaluated using the Dirichlet spectral
data (\ref{Dirichlet_spectral_data}), see the first formula in
(\ref{time_fourier}), it is possible to find, for such $h\,$,the trace,
\bfo
Wh|_{\Sigma} =\sum_{j=1}^\infty \kappa_j^h \phi_j|_{\Sigma}.
\efo
The above considerations give rise to the following lemma.

 \begin{lem}
 \label{Z-lemma}
 Assume $\sigma(-\Delta^D(\overline S)) \cap \sigma(-\Delta^D(M \setminus S))=
 \emptyset$.
 Then the subspaces,
 \bfo
 {\it D}^1_{\tilde \Sigma}:= \{h \in {\it D}^1: Wh|_{\tilde \Sigma}=0\}
 \subset {\it D}^1,
 \efo
 are uniquely determined in terms of the Dirichlet spectral data
 (\ref{Dirichlet_spectral_data}).

Moreover, the wave operator $W$, restricted to ${\it D}^1_{\tilde \Sigma}$,
\bfo
W:{\it D}^1_{\tilde \Sigma} \rightarrow H^1_{\tilde \Sigma},
\quad H^1_{\tilde \Sigma}:= \{a \in H^1(M): a|_{\tilde \Sigma}=0 \},
\efo
is unitary.
 \end{lem}

\subsection{Finding eigenvalues and eigenfunctions in subdomains}
In this subsection we denote by $\tilde S$
one of the manifolds $S_i,\, i=1,2,\,$
$M \setminus S_i,\, i=1,2,\,M \setminus S$ and by $\lambda_n(\tilde S),
\, \phi_n(\cdot; \tilde S)$ we denote the eigenvalues and orthonormal
eigenfunctions of $-\Delta^D(\tilde S)$.
By the
max-min principle,
 \ba
\lambda_n(\tilde S)=\max_{u_1,\dots,u_{n-1}}\min_{u_n}
(\triangledown_gu_n,\triangledown_gu_n)_{L^2(\tilde S)},
 \ea
 where the maximum is taken over
$u_1,u_2,\dots,u_{n-1}\in H^{1}_0(\tilde S)$ and the minimum is taken
over $u_n\in H^{1}_0(\tilde S)$ that satisfies \ba
(u_n,u_p)_{L^2(\tilde S)}=0,\quad p=1,2,\dots,n-1, \quad
(u_p,u_p)_{L^2(\tilde S)}=1,\quad p=1,2,\dots,n.
 \ea
 The minimizer
$u_n(x)$ is then an normalized eigenfunction corresponding to the
eigenvalue $\lambda_n(\tilde S)$.
Now consider the following max-min problem
 \beq
 \label{min_max}
t_n(\tilde \Sigma)=\max_{u_1,\dots,u_{n-1}}\min_{u_n}
(\triangledown_g u_n,\triangledown_g u_n)_{L^2(M)}=\max_{u_1,\dots,u_{n-1}}\min_{u_n}
\sum_{j=1}^\infty \la_j |u_{n,j}|^2,
 \eeq
 where $u_{n,j}$ are the Fourier coefficients of $u_n,$ i.e.
$u_n(x)=\sum_{j=1}^\infty u_{n,j}\phi_j(x)$.
 Here the maximum is taken over
$u_1,u_2,\dots,u_{n-1}\in H^{1}(M)$ satisfying $u_p|_{\tilde \Sigma}=0$ and
the minimum is taken over $u_n\in H^{1}(M)$  satisfuing
$u_n|_{\tilde \Sigma}=0$ and
 \beq
 \label{min_max_1}
& &(u_n,u_p)_{L^2(M)}=\sum_{j=1}^\infty u_{n,j} {\overline u_{p,j}}= 0,
\quad p=1,2,\dots,n-1,
\\ \nonumber
& &(u_p,u_p)_{L^2(M)}= \sum_{j=1}^\infty |u_{p,j}|^2 = 1,\quad p=1,2,\dots,n.
\eeq
Then $t_n(\tilde \Sigma)$ are the eigenvalues of the Dirichlet Laplacian
on the direct sum of $L^2(\tilde S)$ and $L^2(M \setminus \tilde S)$,
so that
\bfo
\{t_n(\tilde \Sigma)\}_{n=1}^\infty= \sigma(-\Delta^D(\tilde S)) \cup
\sigma(-\Delta^D(M \setminus\tilde S)).
\efo
The sequence of the corresponding minimizers,
$u_n(x;\tilde\Sigma)$ consists of orthonormal
eigenfunctions of this operator. However, due to the assumption
$\sigma(-\Delta^D( \tilde S)) \cap
\sigma(-\Delta^D(M \setminus \tilde S))=
 \emptyset$,
any such eigenfunction is equal to $0$ on $M \setminus \tilde S$
 or $\tilde S$.
Thus, any $u_n(x;\tilde\Sigma)$ is either an eigenfunction
of
$-\Delta^D( \tilde S)$ extended by $0$ to $M \setminus \tilde S$, or
an eigenfunction
of
$-\Delta^D(M \setminus  \tilde S)$ extended by $0$ to $\tilde S$.

On the other hand, Lemma \ref{Z-lemma} together with equation
(\ref{time_fourier}) make it possible to evaluate the right-hand
sides in (\ref{min_max}), (\ref{min_max_1}) using the
Dirichlet spectral data. This leads to the following result

\begin{lem}
\label{new}
Let
$\sigma(-\Delta^D( \overline{S_i})),\,  \sigma(-\Delta^D( M \setminus S_i)), i=1,2$
and $\sigma(-\Delta^D(M \setminus S))$ be all disjoint. Then the
Dirichlet spectral data (\ref{Dirichlet_spectral_data}) determine
uniquely the eigenvalues $\la_n(\tilde S),\, n=1, 2, \dots,\,
\tilde S=S, S_i, M\setminus S_i,\, i=1, 2,$ and $ M\setminus S $. They determine also
the the generalized sources $h_n(\tilde S)$ such that
\beq
\label{extended_functions}
W h_n(\tilde S)=
 u^{0,h_n(\tilde S) }(x, 0)=\begin{cases} \phi_n(x;\,\tilde S),\quad x \in \tilde S,
\\
0,\quad x \in M \setminus \tilde S. \end{cases}
 \eeq
 In addition, the
Dirichlet spectral data determine the Fourier coefficients of the
extended eigenfunctions $W h_n(\tilde S)(x)$,
\beq
\label{Fourier}
W h_n(\tilde S)=
 \sum_{j=1}^\infty \kappa_{n, j}(\tilde S) \phi_j(x), \quad x \in M.
\eeq
\end{lem}
We note that these sources $h_n(\tilde S)$ are determined up to a unitary
transformation in the eigenspace corresponding to the eigenvalue $\la_n(\tilde S)$.
\begin{proof}
Recall that, by formula (\ref{time_fourier}) we can evaluate the Fourier
coefficients $u^{0,h}_j(0)$ for any $h \in {\it D}^1$. Thus, we can evaluate
\bfo
\mu_n(\tilde \Sigma) =\max_{h_1,\dots,h_{n-1}}\min_{h_n}\sum_{j=1}^\infty
\la_j |u^{0,h_n}_j(0)|^2,
\efo
where the maximum is taken over
$h_1,h_2,\dots,h_{n-1}\in {\it D}^1_{\tilde \Sigma}$
and the minimum is taken over $h_n\in {\it D}^1_{\tilde \Sigma}$ with
\bfo
\sum_{j=1}^\infty u_{j}^{0,h_n}(0) {\overline u_{j}^{0,h_p}(0)}= 0,
\ p=1,\dots,n-1,\quad
 \sum_{j=1}^\infty |u_{j}^{0,h_p}(0)|^2 = 1,\ p=1,\dots,n.
\efo
It follows from Lemma \ref{Z-lemma} that, for any $n=1,2,\dots,$
$
\mu_n(\tilde \Sigma) =t_n(\tilde \Sigma),
$
providing $\sigma(-\Delta^D(\overline{\tilde S})) \cap \sigma(-\Delta^D(M \setminus \tilde S))=\emptyset$.

Repeating this construction with $\tilde \Sigma$ equal to
$\Sigma$, $\Sigma_1$, and $\Sigma_2$, we obtain the sets
 \ba
\sigma(-\Delta^D(S)) \cup \sigma(-\Delta^D(M \setminus S))=
\sigma(-\Delta^D(S_1)) \cup \sigma(-\Delta^D(S_2)) \cup \sigma(-\Delta^D(M \setminus S)),
 \ea
 and
 \ba
\sigma(-\Delta^D(S_1)) \cup \sigma(-\Delta^D(M \setminus S_1)), \quad \quad
\sigma(-\Delta^D(S_2)) \cup \sigma(-\Delta^D(M \setminus S_2)).
 \ea
As $\sigma(-\Delta^D( \overline{S_i})),\,  \sigma(-\Delta^D( M \setminus S_i)), i=1,2$
and $\sigma(-\Delta^D(M \setminus S))$ are all disjoint, by intersecting the above sets
we find the desired eigenvalues $\la_n(S_i),\, n=1, 2, \dots, i=1, 2$
and $\la_n(M \setminus S_i),\, n=1, 2, \dots, i=1, 2$, as well as
$\la_n(M \setminus S),\, n=1, 2, \dots$.

Then, identifying the corresponding subsequence of $t_n(\tilde \Sigma)$ and
related generalized sources $h_n(\tilde \Sigma)$, we determine, for each $\tilde S$,
the generalized sources $h_n(\tilde S)$ such that
$W h_n(\tilde S)$ are equal to the extended eigenfunctions
(\ref{extended_functions}).

Recalling formula (\ref{time_fourier}), we prove the last part of the Lemma.
\end{proof}

\subsection{Inverse problems in subdomains. Proof of Theorem \ref{thm2}}
Our proof of Theorem  \ref{thm2} is based on Lemma \ref{new}.
Namely, we will show that, having at hand
the eigenvalues $\la_n(S_i),\, i=1,2,$ and $\la_n(M \setminus S)$
and also the Fourier coefficients, $\kappa_{n, j}(S_i), \, i=1,2,\,
\kappa_{n, j}(M \setminus S)$, it is possible to determine, up to an isometry,
the Riemannian manifolds $(S_i, g),\, i=1,2,\,(M \setminus S, g)$. Gluing them
along $\Sigma_i$ we recover $(M, g)$.

Recall that if, for $(N,g),\, \p N \neq \emptyset$, we do know its Dirichlet
eigenvalues $\la_n(N)$ and traces on $\p N$ of the normal derivatives
of the eigenfunctions, $\p_\nu \psi_n|_{\p N}$, then these data determine
$(N,g)$ up to an isometry, see e.g. \cite{KKL}. However, in the case
of the Dirichlet spectral data, we have only the Dirichlet values
$\phi_j|_{\Sigma}$ and, moreover, the convergence of the Fourier series
(\ref{Fourier}) is only in $H^1(M)$ preventing us from identifying
$\p_\nu \phi_n(\tilde S)|_{\tilde \Sigma}$. Therefore, we will use another
approach within the BC-method, described in section 4.1 of \cite{KKL}.
To explain it, consider the intial-boundary value problem in $N  \times \R$,
\begin{equation}
 \label{eq_N}
 \begin{aligned}
 &(\p_t^2-\Delta_g)w_F=0, \quad\hbox{in } N \times \R_+, \\
 &w_F|_{\p N \times \R_+} =F,\quad
  w_F|_{t=0}=0,\, \p_t w_F|_{t=0} =0,
\end{aligned}
 \end{equation}
with $F \in C^\infty_+(\p N \times \R_+)$. The energy, at time $t$,
of the wave $w_F$ is then defined as
\begin{align*}
\label{energy}
 {\it E}(w_F,t)=&\frac{1}{2}
 \int_{N}(|\partial_t w_F(x,t)|^2+
 |\nabla_g
 w_F(x,t)|_g^2)dV_g(x).
 \end{align*}
 It is shown in \cite{KKL} that, given the energy flux
 \bfo
 \Pi(F):=\lim_{t \to \infty}{\it E}(w_F,t),
 \efo
 for any $F \in C^\infty_0(\p N \times \R_+)$,
it is possible to determine $(N, g)$
 up to an isometry.

 Therefore,  Theorem \ref{thm2} is an immediate
corollary of the  following Lemma,
 \begin{lem}
 \label{Energy}
 Let $\Sigma=\Sigma_1 \cup \Sigma_2$ divide $M$ into regions $S_i,\, i=1,2,$ and
 $M \setminus S$. Consider the initial-boundary value problems
  (\ref{eq_N}) with $N$ equal to $\overline{S}_1,\, \overline{S}_2$ and $M \setminus S$.
 Then, assuming that $\sigma(-\Delta^D( \overline{S_i})),\,  \sigma(-\Delta^D( M \setminus S_i)), i=1,2,$
and $\sigma(-\Delta^D(M \setminus S))$ are all disjoint,
the Dirichlet spectral data (\ref{Dirichlet_spectral_data})
uniquely determine
the energy flux $\Pi(F)$  in each of these
subdomains.
\end{lem}
\begin{proof}
We start with the case $N={\overline S_1}$. Then (\ref{eq_N}) takes the form
\begin{equation}
 \label{eq_S}
 \begin{aligned}
 &(\p_t^2-\Delta_g)w_F=0, \quad\hbox{in }  S_1 \times \R_+, \\
 &w_F|_{\Sigma_1 \times \R_+} =F,\quad
  w_F|_{t=0}=0,\, \p_t w_F|_{t=0} =0,
\end{aligned}
 \end{equation}
  with $F \in C^\infty_+(\Sigma_1 \times \R_+)$.
 Let us first show that, for any such
$F$,
there exists a unique $h_F \in C^\infty_+(\Sigma_1 \times \R_+)$
 such that
 $
 w_F=u^{0,h}|_{S_1},
 $
 where
 \beq
 \label{h}
 h=\begin{cases} h_F,\quad x \in \Sigma_1,
\\
0,\quad x \in \Sigma_2. \end{cases}
 \eeq
To this end we consider, in addition to (\ref{eq_S}), the problem
\bfo
 \begin{aligned}
 &(\p_t^2-\Delta_g)w_F^c=0, \quad\hbox{in } (M \setminus S_1)\times \R_+, \\
 &w_F^c|_{\Sigma_1 \times \R_+} =F,\quad
  w_F^c|_{t=0}=0,\, \p_t w_F^c|_{t=0} =0,
\end{aligned}
\efo
and introduce the function $u$,
\beq \label{wF}
u=\begin{cases} w_F,\quad x \in S_1 \times \R_+,
\\
w_F^c,\quad x \in (M \setminus S_1)\times \R_+. \end{cases}
 \eeq
Then $u$ solves the transmission problem (\ref{eq_transmission_t})
with $f=0$ and
 \bfo
 h=\begin{cases} h_F=
\left(\p_\nu w_F^c-  \p_\nu w_F  \right)|_{\Sigma_1 \times \R_+},
\quad (x,t) \in \Sigma_1 \times \R_+
\\
0,\quad (x,t) \in
\Sigma_2
\times \R_+,
\end{cases}
 \efo
 i.e. $u=u^{0, h}$.
Using considerations similar to those in the proof of Lemma
\ref{lem_neumann-to-dirichlet}, we show the uniqueness of such $h$.

By Corollary \ref{trace}, we can then find, for
$h \in C^\infty_+(\Sigma_1 \times \R_+)$,
%Observe now that, using formula (\ref{time_fourier}), we can evaluate
%the Fourier coefficients, $u^{0,h}_j(t)$ of  $u^{0,h}(t),\, t>0$.
%Then, using the Dirichlet spectral data (\ref{Dirichlet_spectral_data}),
%we can find,
% due to the
%convergence of the eigenfunction decomposition in $H^1(M)$,
\beq
\label{nonstationary}
\Lambda_1 h = u^{0,h}|_{\Sigma_1 \times \R_+}.
\eeq
As shown earlier, for $h=0$ on $\Sigma_2$, the operator $\Lambda_1$ is an invertible
operator in $C^\infty_+(\Sigma_1 \times \R_+)$. Thus, we can use
equation (\ref{nonstationary}) with the right-hand side being $F$,
to uniquely determine $h_F$. Observe that the extended eigenfunctions
$\{W \left(h_n(S_1)\right)\}_{n=1}^\infty$ together with
$\{W \left(h_k(M \setminus S_1)\right)\}_{k=1}^\infty$ form an orthonormal
basis in $L^2(M)$. Thus, using (\ref{Fourier}),  we can evaluate
the Fourier coefficients of $u^{0,h}(\cdot, t)$
with respect to this basis,
\bfo
u^{0,h}(x, t)= \sum_{n=1}^\infty w_{n, F}(t)\, W \left(h_n(S_1)\right)(x)+
\sum_{k=1}^\infty w_{k, F}^c(t) \, W \left(h_k(M \setminus S_1)\right)(x),
\efo
where the index $F$ indicates that $h$ is of form (\ref{h}).
This expansion, together with the definition of $u$,
see (\ref{wF}), shows that
\beq \label{neww} \nonumber
& & w_F(x, t) = \Theta_1(x)u^{0,h}(x, t)=  \sum_{n=1}^\infty w_{n, F}(t)
\, W \left(h_n(S_1)\right)(x),
\\
& &\|w_F(t)\|^2_{L^2(S_1)}=\sum_{n=1}^\infty |w_{n, F}(t)|^2,
\eeq
where $\Theta_1$ is the characteristic function of $S_1$ and we identify
functions in $S_1$ with their extensions by $0$ to $M \setminus S_1$.

Similarly,
\bfo
\p_tw_F(x, t)=\sum_{n=1}^\infty \p_tw_{n, F}(t)\,W \left(h_n(S_1)\right)(x),
\quad
\|\p_tw_F(t)\|^2_{L^2(S_1)}=\sum_{n=1}^\infty |\p_tw_{n, F}(t)|^2
\efo
The last step of the construction is based on the observation that,
when $F \in C^\infty_0(\Sigma_1 \times \R_+)$, then
$w_F \in H^1_0(S_1)$
for large $t$ and, therefore, the eigenfunction
expansions (\ref{neww}) converge in $H^1(S_1)$. Thus, we can find
\bfo
\|\nabla_g w_F(\cdot, t)\|_{L^2(S_1)}^2=
\sum_{n=1}^\infty \la_n(S_1) |w_{n, F}(t)|^2.
\efo
Combining the above two equations, we find the energy flux $\Pi(F)$.

Similar considerations, with $\overline{S}_2$ and $M \setminus S_2$ and $M \setminus S$
and $\overline{S}$ instead of $\overline{S}_1$ and $M \setminus S_1$, show the possibility to
evaluate the energy flux when $N=\overline{S}_2$ and $M \setminus S$.
\end{proof}

\section{Some remarks on further
generalizations and open problems}

\label{sec_measurements}

 \begin{rem} The inverse problem with the
Cauchy or Dirichlet spectral data on a compact
Riemannian manifold without boundary studied in Sections
\ref{sec_cauchy},
\ref{sec_dirichlet} can be generalized to the problem when
 $M$ has boundary $\pM \neq \emptyset$ with e.g. Neumann (or Dirichlet)
boundary condition. Then
(\ref{Cauchy_spectral_data}) and (\ref{Dirichlet_spectral_data})
consists of the Neumann eigenvalues
of $M$ and the Cauchy or Dirichlet data of the eigenfuctions on a
closed hypersurface $\Sigma\subset M$, $\p M\cap \Sigma=\emptyset$.
The methods to solve these problems are very similar to the described
and are based on the version of the BC-method with data on a part
of boundary, see section 4.4 \cite{KKL}.

 \end{rem}

\begin{rem}

Consider the non-stationary Green function, $G(x,y;t)$, for
the acoustic wave equation in
$(M, g)$, where $M$ may have a non-trivial boundary,
 \begin{align*}
 &(\p_t^2-\Delta_g)G(x,y;t)=\delta_y(x)\delta(t),\quad \text{in } M\times\R,\\
 &G(x,y;t)|_{\p M \times\R_+}=0,\quad
 G(x,y;t)|_{t<0}=0,
 \end{align*}
where $y\in\Sigma$, $\Sigma\subset M$ such that
$\Sigma\cap\p M=\emptyset$ and the boundary condition is void
when $\p M= \emptyset$. It often happens in practice,
for example in geophysics or ultrasound imaging,
that one can measure $G(x,y;t)$ for $x$ again running only over $\Sigma$.
Thus, we come to the inverse problem with non-stationary data
being $G(x,y;t),\, x, y \in \Sigma,\, t>0$.

Taking the Fourier transform of the given $ G(x,y;t)$  in
 $t$, we find the Green function $G(x,y;k)$, cf. \cite{KKLM},
 \[
 (-k^2-\Delta_g)G(x,y;k)=\delta_y(x), \quad G(x,y;t)|_{\p M}=0.
 \]
Thus, from practical measurements we can find
$G(x,y;k)$ for $x,y \in \Sigma,\, k \in \R$.

Note that $G(x,y;k),\, x, y \in \Sigma,$ is the integral kernel of
 a meromorphic, with respect to $k \in \C$, operator-valued function in
$L^2(\Sigma)$. In terms of the eigenvalues and eigenfunctions
of $-\Delta_g$,
 \[
 G(x,y;k)|_{\Sigma \times \Sigma}=\sum_{j=1}^\infty\frac{1}{k^2-\lambda_j}\phi_j(x)|_\Sigma\phi_j(y)|_\Sigma,
 \]
 where the convergence in the right-hand side is understood in the sense of operators
 in $L^2(\Sigma)$.
So finding the poles of $G(x,y;k)$ is equivalent to determination
of $\sqrt{\lambda_j}$. At the pole $\sqrt{\lambda_j}$, the residue
is given by the integral operator with the kernel,
 \[
\text{Res}(G(x,y;\cdot),\sqrt{\lambda_j})=\frac{1}{2\sqrt{\lambda_j}}
\sum_{l:\lambda_l=\lambda_j}\phi_l(x)|_\Sigma\phi_l(y)|_\Sigma.
 \]
The knowledge of this kernel  allows us to find the
functions $\phi_j|_{\Sigma}$   up to an orthogonal transformation
in the eigenspace corresponding to $\la_j$, see \cite{KKLM}, \cite{Rom}.

Therefore, the dynamic inverse data $G(x,y;t),\, x, y \in \Sigma,\, t>0$
makes it possible to find the Dirichlet spectral data on $\Sigma$.

\end{rem}

\begin{rem}
 As shown in section \ref{sec_dirichlet}, the Dirichlet data
 (\ref{Dirichlet_spectral_data}) determine $(M, g)$ when the spectra
$ \sigma(-\Delta^D( \overline{S_i})),\,  \sigma(-\Delta^D( M \setminus S_i)), i=1,2,$
and $\sigma(-\Delta^D(M \setminus S))$ are all disjoint.
It is interesting to understand whether this condition can be removed.

In general, it is important to find if the Dirichlet spectral data
 (\ref{Dirichlet_spectral_data}) determine $(M, g)$  even in the case when $\Sigma$ is connected.
\end{rem}

\section{Acknowledgements}

%\HOX{Please add your projects}

The research of K.K. was financially supported by the
Academy of Finland (project 108394) and the research of M.L. by the Academy of Finland Center of Excellence programme 213476.

 \bibliographystyle{amsplain}
 \bibliography{closed_manifold_2007_06m_10d}

\end{document}